\documentclass[11pt,twoside]{amsart}

\usepackage{amsfonts}

\title[K\"ulshammer-like Hochschild homology invariants]{Hochschild
homology invariants of K\"ulshammer type of derived categories}
\author{Alexander Zimmermann}
\address{Universit\'e de Picardie,\newline
D\'epartement de Math\'ematiques et
\newline LAMFA (UMR 6140 du CNRS),
\newline 33 rue St Leu,
\newline 80039 Amiens Cedex 1,
\newline France}
\email{alexander.zimmermann@u-picardie.fr}
\urladdr{http://www.lamfa.u-picardie.fr/alex/azim.html}
\date{January 6, 2010; revised May 10, 2010}
\thanks{I gratefully acknowledge financial support from a grant
``partenariat Hubert Curien PROCOPE''}
\dedicatory{Herrn Klaus W. Roggenkamp zum 70. Geburtstag}



\newtheorem{Theo1}{{Theorem}}
\newtheorem*{Theo2}{{Theorem}}
\newtheorem{Lemma1}{{Lemma}}
\newtheorem{Def1}{{Definition}}
\newtheorem{Prop1}[Lemma1]{{Proposition}}
\newtheorem{Claim1}[Lemma1]{{Claim}}
\newtheorem{Rem1}{{Remark}}[section]
\newtheorem{Cor1}[Lemma1]{{Corollary}}
\newtheorem{Ex1}{{Example}}

\newenvironment{Lemma}{\begin{Lemma1}}{\end{Lemma1}}
\newenvironment{Def}{\begin{Def1}\em}{\end{Def1}}
\newenvironment{Prop}{\begin{Prop1}}{\end{Prop1}}
\newenvironment{Rem}{\begin{Rem1}\rm}{\end{Rem1}}
\newenvironment{Theorem}{\begin{Theo1}}{\end{Theo1}}

\newenvironment{Remark}{\begin{Rem1}\em}{\end{Rem1}}

\newcommand{\uar}{\uparrow}
\newcommand{\dar}{\downarrow}
\newcommand{\lra}{\longrightarrow}
\newcommand{\lla}{\longleftarrow}
\newcommand{\ra}{\rightarrow}
\newcommand{\sdp}{\times\kern-.2em\vrule height1.1ex depth-.05ex}
\newcommand{\epi}{\lra \kern-.8em\ra}

\newcommand{\N}{{\mathbb N}}

\newcommand{\T}{{\mathbb T}}
\newcommand{\B}{{\mathbb B}}

\newcommand{\dickebox}{{\vrule height5pt width5pt depth0pt}}

 \normalbaselines
\begin{document}

\begin{abstract} For a perfect field $k$ of characteristic $p>0$
and for a finite dimensional symmetric $k$-algebra $A$ K\"ulshammer studied a
sequence of ideals of the centre of $A$ using the $p$-power map on degree
$0$ Hochschild homology. In joint work with Bessenrodt and Holm we removed
the condition to be symmetric by passing through the trivial extension algebra.
If $A$ is symmetric then the dual to the K\"ulshammer ideal structure was
generalised to higher Hochschild homology in earlier work \cite{gerstenhaber}.
In the present paper we follow this program and
propose an analogue of the dual to the K\"ulshammer ideal structure
on the degree $m$ Hochschild homology
theory also to not necessarily symmetric algebras.
\end{abstract}

\maketitle

\section*{Introduction}

Let $k$ be a perfect field of characteristic $p>0$ and let $A$ be a
finite dimensional $k$-algebra.
For a symmetric algebra $A$ K\"ulshammer introduced in \cite{Kueldeutsch}
a descending sequence of ideals $T_n(A)^\perp$ of the centre of $A$ satisfying
various interesting properties. In \cite{HHKM} H\'ethelyi, Horv\'ath,
K\"ulshammer and Murray continued to study this sequence of ideals,
proved their invariance under Morita equivalence, and showed further
properties linked amongst others to the Higman ideal of the algebra.
In order to show these properties K\"ulshammer introduced a sequence
of mappings $\zeta_n:Z(A)\lra Z(A)$
with image being $T_n(A)^\perp$. In a completely dual procedure he
introduced mappings $\kappa_n:A/KA\lra A/KA$
again in case $A$ is symmetric. These mappings encode many representation
theoretic informations, in particular in case $A$ is a group algebra. The reader
may refer to K\"ulshammer's original articles \cite{Kueldeutsch}
and \cite{Kuelprog} for more  details. Moreover, in a recent
development it could be proved in \cite{kuelsquest} that these mappings
actually are an invariant of the derived category of $A$.
This fact could be used to
distinguish derived categories in very subtle situations, such as some
parameter questions for blocks of dihedral or semidihedral type in
joint work with Holm \cite{hz-tame}, such as very delicate questions
to distinguish two families
of symmetric algebras of tame domestic representation type by Holm and
Skowro\'nski \cite{Holmskow}, or in joint work with Holm \cite{typeL}
to distinguish derived
equivalence classes of deformed
preprojective algebras of generalised  Dynkin type, as defined by
Bia\l kowski, Erdmann and Skowro\'nski \cite{BES} (modified slightly for type $E$;
cf \cite[page 238]{periodic}) with respect to different
parameters.
Most recently in joint work \cite{stableHochschild}
with Yuming Liu and Guodong Zhou we
gave a version which is invariant under stable equivalences of Morita type.
This approach gives a link to Auslander's conjecture on the invariance of
the number of simple modules under stable equivalences of Morita type.
K\"onig, Liu and Zhou \cite{KLZ} continued the work with focus being higher Hochschild
and cyclic homology invariants.

In joint work with Bessenrodt and Holm \cite{nonsymmetric} we could get
rid of the assumption that $A$ needs to be symmetric for the definition
of the images of $\zeta_n$. We used the trivial extension algebra $\T A$
of $A$, computed the image of $\zeta_n$ for $\T A$ and interpreted the result
purely in terms of $A$. The derived invariant K\"ulshammer's ideal structure
of the centre of $A$ becomes available for any finite dimensional algebra over
perfect fields of positive characteristic.

Since the centre of the algebra is the degree $0$ Hochschild cohomology and
$A/KA$ is the $0$-degree Hochschild homology, it is natural to try to
understand the K\"ulshammer construction for higher degree Hochschild (co-)homology.
Denote by $HH^m(A)$ the $m$-degree
Hochschild cohomology of $A$ and by $HH_m(A)$ the $m$-degree Hochschild homology
of $A$.
We defined in \cite{gerstenhaber} for a symmetric algebra $A$ certain mappings
$\kappa_n^{(m);A}:HH_{p^nm}(A)\lra HH_m(A)$
in completely analogous manner as K\"ulshammer's $\kappa_n$. Moreover
we showed that $\kappa_n^{(0);A}=\kappa_n$, and proved its derived invariance as well
(cf Theorem~\ref{derivedkappacup} below).

In Definition~\ref{Definitionofkappahat} of Section~\ref{Kuelshammerfornonsymmetric} we propose a mapping
$\hat \kappa_n^{(m);A}$, denoted also $\hat \kappa_n^{(m)}$
if no confusion may occur, in analogy of $\kappa_n^{(m);A}$ for algebras which are not necessarily symmetric.
Further, we prove that $\hat \kappa_n^{(m)}$ is invariant under a derived equivalence:
If
$F:D^b(A)\simeq D^b(B)$ is a standard equivalence
then it induces coherent isomorphisms of Hochschild homology
groups $HH_m(F):HH_m(A)\lra HH_m(B)$
and we get
$HH_m(F)\circ \hat\kappa_n^{(m);A}=\hat\kappa_n^{(m);B}\circ HH_{p^nm}(F)$ for all
positive integers $m$ and~$n$.

Finally we give an elementary example in order to illustrate that these
invariants are computable and to show that they are not trivial.

The paper is organised as follows. In Section~\ref{Hochschilddefinitions} we
recall some known facts about Hochschild (co-)homology.
In Section~\ref{Kuelshammerrappel} the known constructions concerning
K\"ulshammer ideals in degree $0$ for symmetric and non
symmetric algebras, using trivial extension algebras are recalled, as
well as the $m$-degree Hochschild cohomology generalisation
for symmetric algebras.
Section~\ref{barcomplexsection} recalls that Hochschild homology is functorial,
a statement which is used at a prominent point in our construction.
In Section~\ref{Kuelshammerfornonsymmetric} finally
we provide a K\"ulshammer mapping for not necessarily symmetric algebras
and higher Hochschild homology.
We show in
Theorem~\ref{derivedinvarianceofkappahat} that the new
mapping $\hat\kappa_n^{(m)}$ is an invariant of the derived category.
Section~\ref{dualnumbers} is devoted to a detailed computation of an example.

\bigskip

{\bf Acknowledgement:} Part of this work was done in October
2007 during a visit at Beijing Normal University and Nanjing Normal University.
I wish to express my gratitude to
Changchang Xi and to Jiaqun Wei for their invitation and their kind hospitality.

\section{Setup of the basic tools}

\subsection{Hochschild homology and cohomology}
\label{Hochschilddefinitions}

Recall some well-known facts from Hoch\-schild (co-)homology. We refer to
Loday \cite{Loday} for a complete presentation of this theory.
Let $A$ be a finitely generated
$k$-algebra for a commutative ring $k$, and suppose $A$ is
projective as $k$-module.
Then for every $A\otimes_kA^{op}$-module $M$ we define
the Hochschild homology of $A$ with values in $M$ as
$$HH_n(A,M):=Tor_n^{A\otimes_kA^{op}}(A,M)$$
and Hochschild cohomology of $A$ with values in $M$ as
$$HH^n(A,M):=Ext_{A\otimes_kA^{op}}^n(A,M).$$
Often and frequently throughout the paper, for a $k$-algebra $R$ we abbreviate
$R^e:=R\otimes_kR^{op}$.
There is a standard way to compute Hochschild (co-)homology
by the bar resolution. We recall its definition in order to set
up the notations.

The bar complex $\B A$ is given by
$$(\B A)_n:=A^{(\otimes_k)^{n+1}}:=
\underbrace{A\otimes_kA\otimes_k\dots\otimes_kA}_{n+1\mbox{ \scriptsize
factors}}$$
and differential
$$d(a_0\otimes\dots a_n)=
\sum_{j=0}^{n-1}(-1)^ja_0\otimes\dots\otimes
a_{j-1}\otimes a_ja_{j+1}\otimes a_{j+2}\otimes \dots\otimes a_n.$$
One observes that $d^2=0$, that $(\B A)_n$ is a free $A\otimes_kA^{op}$-module
for every $n> 0$ and that the complex $(\B A,d)$ is
a free resolution of $A$ as a $A\otimes_kA^{op}$-module, called the bar resolution.
Hence, $$HH_n(A,M)\simeq H_n(\B A\otimes_{A\otimes A^{op}}M)\mbox{ and }
HH^n(A,M)\simeq H^n(Hom_{A\otimes A^{op}}(\B A,M)).$$

Moreover, the bar complex is functorial in the sense that whenever
$f:A\lra B$ is a homomorphism of $k$-algebras, then
$$(\B A)_n\ni a_0\otimes\dots \otimes a_n\mapsto f(a_0)\otimes\dots \otimes
f(a_n)\in \B B$$
induces a morphism of complexes $\B A\lra \B B$.

\subsection{The K\"ulshammer ideals; constructions and their Hochschild
generalisation}

\label{Kuelshammerrappel}

We recall the constructions introduced by K\"ulshammer in \cite{Kueldeutsch}
and \cite{Kuelprog} as well as their generalisation introduced in \cite{gerstenhaber}
to higher degree Hochschild cohomology and their generalisations to $0$-degree
Hochschild cohomology for not necessarily symmetric algebras.

\subsubsection{The K\"ulshammer ideal theory for algebras which are symmetric}
\label{Kuelshammeridealsknownfacts}

Let $k$ be a perfect field of characteristic $p>0$
and let $A$ be a $k$-algebra. Then
$$KA:=<ab-ba\;|\;a,b\in A>_{k-\mbox{\scriptsize space}}$$
is a $Z(A)$-module and for all $n\in \N$
$$T_n(A):=\{x\in A\;|\;x^{p^n}\in KA\}$$
is a $Z(A)$ submodule of $A$. If $A$ is symmetric with symmetrising form
$$<\;,\;>:A\otimes_kA\lra k$$ taking orthogonal spaces, since $KA^\perp=Z(A)$,
$T_n(A)^\perp$ is an ideal of $Z(A)$.

Given an equivalence of triangulated
categories of standard type (cf Section~\ref{derivedequivonhomology} below)
$D^b(A)\lra D^b(B)$ between two symmetric
$k$-algebras $A$ and $B$, Rickard showed in \cite{Ri3} that
this equivalence induces an isomorphism $Z(A)\simeq Z(B)$
and if $k$ is a perfect field, it is shown in \cite{kuelsquest} that
$T_n(A)^\perp$ is mapped by this isomorphism to $T_n(B)^\perp$.

Moreover $A$ is symmetric if and only if $A\simeq A^*:=Hom_k(A,k)$ as
$A\otimes_kA^{op}$-modules.
Since we supposed $k$ to be perfect the Frobenius map is invertible and let
$k\ni \lambda\mapsto\lambda^{p^{-1}}\in k$ be its inverse.
Denote the $n$-th iterate by $p^{-n}:=\left(p^{-1}\right)^n$.
Since $A/KA\ni b\mapsto b^p\in A/KA$ is a well-defined additive mapping,
semilinear with respect to the Frobenius mapping, one defines an
element in $Hom_k(A,k)$ by
$$a\mapsto \left(b\mapsto <a,b^{p^n}>^{p^{-n}}\right)$$
and so for every $a\in Z(A)$ there is a unique element $\zeta_n(a)\in A$
so that $$<a,b^{p^n}>=<\zeta_n(a),b>^{p^n}\forall\;b\in A.$$

Similarly, using that $Z(A)^\perp=KA$, one gets a mapping
$\kappa_n:A/KA\lra A/KA$ satisfying
$$<a^{p^n},b>=<a,\kappa_n(b)>^{p^n}\forall\;a\in Z(A), b\in A/KA.$$

In \cite{gerstenhaber} this construction was generalised to a structure on
Hochschild homology which is invariant under equivalences of derived categories.
More precisely, let $A$ be a {\em symmetric} finite dimensional
$k$-algebra for $k$
a perfect field of characteristic $p>0$.
Then, denoting by $\B A$ the bar resolution of $A$ and $A^e=A\otimes_kA^{op}$,
$$Hom_k(\B A\otimes_{A^e}A,k)
\simeq Hom_{A^e}(\B A,Hom_k(A,k))\simeq Hom_{A^e}(\B A,A)$$
and so $$Hom_k(HH_m(A,A),k)\simeq HH^m(A,A)$$
which induces a non degenerate pairing
$$<\;,\;>_m:HH^m(A,A)\times HH_m(A,A)\lra k\;.$$
By a construction analogous to the one for the pairing
$$<\;,\;>:Z(A)\times A/KA\lra k$$
one gets a mapping
$$\kappa_n^{(m);A}:HH_{p^nm}(A,A)\lra HH_m(A,A)$$
denoted simply by $\kappa_n^{(m)}$ if no confusion may occur, satisfying
$$<z^{p^n},x>_{p^nm}=\left(<z,\kappa_n^{(m)}(x)>_{m}\right)^{p^n}$$
for all $x\in HH_{p^nm}$ and $z\in HH^m(A,A)$.

\subsubsection{On derived equivalences and Hochschild (co-)homology}

\label{derivedequivonhomology}

In \cite{gerstenhaber} we studied invariance properties of these
mappings with respect to equivalences between derived categories.
Rickard showed in \cite{Ri3} that whenever $A$ and $B$ are $k$-algebras over a
field $k$, and whenever $D^b(A)$ and $D^b(B)$ are equivalent as triangulated categories, then there is a complex $X\in D^b(B\otimes_kA^{op})$ so that
$X\otimes_A^{\mathbb L}-:D^b(A)\lra D^b(B)$ is an equivalence. Such an equivalence is called of standard type and $X$ is called a two-sided tilting complex. It is shown in \cite{Ri3} that every equivalence of triangulated categories $D^b(A)\lra D^b(B)$ coincides with a suitable standard one on objects.

A tilting complex $T$ is a complex in $K^b(A)$, the homotopy category of bounded complexes of projective $A$-modules, satisfying  $Hom_{D^b(A)}(T,T[i])=0$ for all $i\neq 0$, and $add(T)$ generates $K^b(A)$ as triangulated category.
Two $k$-algebras $A$ and $B$ have equivalent derived categories $D^b(A)$ and $D^b(B)$ if and only if there is a tilting complex $T$ in $D^b(A)$ with endomorphism ring $B$. Given an equivalence $D^b(B)\lra D^b(A)$, the image of the regular $B$-module is a tilting complex.
Hence, by the above, for every tilting complex $T$ in $D^b(A)$ there is a two-sided tilting complex $X$ in $D^b(A\otimes_kB^{op})$ so that $X\simeq T$ in $D^b(A)$.

Let $F_X=X\otimes_A^{\mathbb L}-:D^b(A)\lra D^b(B)$ be an equivalence
of derived categories with quasi-inverse $Y\otimes^{\mathbb L}-$ (cf
Rickard \cite{Ri3}). If $A$ is symmetric, then
$B$ is symmetric as well (cf \cite{rogquest}
for general base rings, and \cite{Ri3} for fields) and
$$X\otimes_A^{\mathbb L}-\otimes_A^{\mathbb L}Y:
D^b(A\otimes_kA^{op})\lra D^b(B\otimes_kB^{op})$$
is an equivalence again. This equivalence induces an isomorphism
$$HH_*(F_X):HH_*(A)\lra HH_*(B)$$
made explicit in \cite{gerstenhaber}.
Recall the precise mapping, defined on the level of complexes, which will be needed later.
\begin{eqnarray*}
\B A\otimes_{A\otimes_kA^{op}}(X\otimes_BY)&\lra&(Y\otimes_A\B A\otimes_AX)\otimes_{B\otimes_kB^{op}}B\\
u\otimes(x\otimes y)&\mapsto&(y\otimes u\otimes x)\otimes 1
\end{eqnarray*}
Then $HH_*(F_X)$ is the mapping induced on the homology.

We recall
one of the main results of \cite{gerstenhaber}.
This result will be one of our main
ingredients in the proof of Theorem~\ref{derivedinvarianceofkappahat}.

\begin{Theorem}\label{derivedkappacup}\cite{gerstenhaber}
Let $A$ be a finite dimensional symmetric $k$-algebra
over the perfect field $k$ of
characteristic $p>0$. Let $B$ be a second algebra such that
$D^b(A)\simeq D^b(B)$ as triangulated categories,
and let $m\in\N$.
Then, there is a standard equivalence $F:D^b(A)\simeq D^b(B)$,
and any such standard equivalence induces an isomorphism
$HH_m(F):HH_m(A,A)\lra HH_m(B,B)$ of all Hochschild homology
groups, satisfying
$$HH_{m}(F)\circ\kappa_n^{(m),A}=
\kappa_n^{(m),B}\circ HH_{p^nm}(F)\;.$$
\end{Theorem}

We take the opportunity to mention that in \cite{gerstenhaber} in the statement of
the theorem the condition
on the field $k$ to be perfect is unfortunately missing.
However, this was a general assumption throughout \cite{gerstenhaber}.

\subsubsection{Trivial extension of an algebra; the K\"ulshammer ideal
theory in the general case}
\label{TA}

Let $k$ be a commutative ring and let $A$ be a $k$-algebra.
Then $Hom_k(A,k)$ is an $A\otimes_kA^{op}$-module by the action
$$\left((a,b)\cdot f\right)(x):=
(a\cdot f\cdot b)(x):=f(bxa)\;\;\forall a,x\in A;b\in A^{op};f\in Hom_k(A,k).$$

Now, the $k$-vector space $$\T A:= A\oplus Hom_k(A,k)$$ becomes
a $k$-algebra by putting
$$(a,f)\cdot (b,g):= (ab,ag+fb)\;\;\forall a,b\in A;f,g\in Hom_k(A,k)\;.$$
Then
\begin{eqnarray*}
A&\stackrel{\iota_A}{\lra}&\T A\\
a&\mapsto&(a,0)
\end{eqnarray*}
is a ring homomorphism. Moreover $\{0\}\oplus Hom_k(A,k)$ is a two-sided ideal
of $\T A$ and the canonical projection
\begin{eqnarray*}
\T A&\stackrel{\pi_A}{\lra}&A\\
(a,f)&\mapsto&a
\end{eqnarray*}
is a splitting for $\iota_A$; i.e. $\pi_A\circ\iota_A=id_A$.
The algebra $\T A$ is always symmetric via the bilinear form
$$\T A\otimes_k\T A\ni <(a,f);(b,g)>\mapsto f(b)+g(a)\in k$$
which induces, for $k$ a field and $A$ finite dimensional over $k$, an
isomorphism
\begin{eqnarray*}
\T A&\lra&Hom_k(\T A,k)\\
(a,f)&\mapsto&\left((b,g)\mapsto f(b)+g(a)\right)
\end{eqnarray*}
of $\T A\otimes_k\T A^{op}$-modules.

In the joint paper \cite{nonsymmetric} with Bessenrodt and Holm
it is shown that the K\"ulshammer ideals satisfy
$$T_n(\T A)^\perp=\{0\}\oplus Ann_{Hom_k(A,k)}(T_n(A))$$
for all $n\geq 1$, and therefore for all $n\geq 1$
the $Z(A)$-modules $Ann_{Hom_k(A,k)}(T_n(A))$ are invariants
under derived equivalences.

\begin{Remark} \label{findbackoriginal}
It should be noted that for algebras $A$ which are already symmetric,
via the symmetrising form $<\;,\;>$ we get an isomorphism of
vector spaces $\lambda:Z(A)\lra Hom_k(A/KA,k)$. Now, every linear form on $A$
equals a form of the shape $<z,->$ for some $z\in Z(A)$ and so
$$\lambda^{-1}(Ann_{Hom_k(A,k)}(T_n(A)))=
\{z\in Z(A)\;|\;<z,T_n(A)>=0\}=T_n(A)^\perp.$$
We find back our original result.
\end{Remark}

\subsection{Relating Hochschild homology of an algebra and of its trivial extension}
\label{barcomplexsection}

Let $R$ and $S$ be $k$-algebras and let $\alpha:R\lra S$ and $\beta:S\lra R$ be
algebra homomorphisms. Then it is well-known
(cf e.g. Loday \cite[Chapter 1, Section 1.1.4]{Loday})
that Hochschild homology is functorial, i.e. $\alpha$ and $\beta$
induce mappings
$$HH_*(\alpha):HH_*(R)\lra HH_*(S)\mbox{ and }HH_*(\beta):HH_*(S)\lra HH_*(R)$$
so that $$HH_*(id_R)=id_{HH_*(R)}\mbox{
as well as
}HH_*(\beta)\circ HH_*(\alpha)=HH_*(\beta\circ\alpha).$$
In particular in case of the mappings $\iota_A:A\lra \T A$ and $\pi_A:\T A\lra A$
from an algebra $A$ to its trivial extension $\T A$ we get induced mappings
giving a split projection
$$HH_n(\pi_A):HH_n(\T A)\lra HH_n(A)$$
and a split injection
$$HH_n(\iota_A):HH_n(A)\lra HH_n(\T A).$$

This can be defined on the level of complexes.
As seen in Section~\ref{Hochschilddefinitions}
the algebra homomorphism $\alpha: R\lra S$ induces a morphism of complexes
\begin{eqnarray*}
\B\alpha:\B R&\lra& \B S\\
(x_0\otimes\dots\otimes x_n)&\mapsto&
(\alpha(x_0)\otimes\dots\otimes \alpha(x_n))
\end{eqnarray*}
and likewise for $\beta:S\lra R$.
Also on the complex computing Hochschild homology an analogous
morphism of complexes is defined by
\begin{eqnarray*}
\B R\otimes_{R^e}R&\lra& \B S\otimes_{S^e}S\\
(x_0\otimes\dots\otimes x_n)\otimes x_{n+1}&\mapsto&
(\alpha(x_0)\otimes\dots\otimes \alpha(x_n))\otimes \alpha(x_{n+1})
\end{eqnarray*}
and likewise for $\beta.$
Its homology induces a mapping $H_*\B\alpha:HH_*(R)\lra HH_*(S)$
which is easily seen to be $HH_*(\alpha)$ obtained above.

%


In particular in our situation
$$H_*\B\pi_A\circ H_*\B\iota_A=id_{\B A\otimes_{A^e}A}$$
and it is clear that
$$H_*\B\pi_A=HH_*(\pi_A)\mbox{ as well as }H_*\B\iota_A=HH_*(\iota_A).$$

\medskip

We reproved the following proposition which was obtained in a much larger context and generality, and much more sophisticated methods, by Cibils, Marcos, Redondo and Solotar in \cite[Theorem 5.8]{CMRS}.

\begin{Prop}\label{directfactor} \cite[Theorem 5.8]{CMRS}
Let $k$ be a field and let $A$ be a $k$-algebra. Then the canonical embedding
$A\lra \T A$ induces a canonical embedding of
$HH_*(A)$ as a direct factor of $HH_*(\T A)$.
\end{Prop}

\section{K\"ulshammer-like Hochschild invariants for non symmetric algebras}
\label{Kuelshammerfornonsymmetric}

Let $k$ be a perfect field of characteristic $p>0$ and let $A$ be a
$k$-algebra. Recall from Section~\ref{Kuelshammeridealsknownfacts} that in case
$A$ is symmetric, we defined in \cite{gerstenhaber} mappings
$$\kappa_n^{(m)}:HH_{p^nm}(A)\lra HH_m(A)$$
satisfying
$$\left(<x,\kappa_n^{(m)}(y)>_{m}\right)^{p^n}=<x^{p^m},y>_{p^nm}\;
\forall \;x\in HH^{m}(A);y\in HH_{p^nm}(A).$$
We should mention that the pairing $<\;,\;>_m$ is defined only on the
(co-)homology of finite dimensional symmetric algebras by the isomorphism
$$
\B A\otimes_{A^e}A\simeq Hom_k(Hom_k(\B A\otimes_{A^e}A,k),k)\simeq
Hom_k(Hom_{A^e}(\B A,A),k)$$
whose homology gives an isomorphism
$$HH_*(A)\simeq Hom_k(HH^*(A),k)\;.$$
Observe that the double dual is the identity for finite dimensional vector
spaces only, and hence here we use the fact that finite dimensional algebras
have projective bimodule resolutions with finite dimensional homogeneous
components.

\begin{Def}\label{Definitionofkappahat}
Let $A$ be a finite dimensional (not necessarily symmetric) $k$-algebra.
Then put
$$\hat\kappa_n^{(m);A}:=
HH_m(\pi_A)\circ \kappa_n^{(m);\T A}\circ HH_{p^nm}(\iota_A)$$
\end{Def}

\bigskip

We shall prove now that the invariant $\hat\kappa_n^{(m)}$ is an
invariant of the derived category of $A$ in the same sense as
it was proved for symmetric algebras (cf
Theorem~\ref{derivedkappacup}).

Rickard showed in \cite{Rickard} that an equivalence between the derived categories of two $k$-algebras induces an equivalence between the derived categories of their trivial extension algebras. We need the following improvement of his result, which seems to be of interest in its own right.

\begin{Prop}\label{twosidedtrivialextension}
Let $A$ and $B$ be finite dimensional $k$ algebras and let $T$ be a tilting complex in $D^b(A)$ with endomorphism ring $B$.
\begin{enumerate}
\item
(Rickard \cite{Rickard}) Then $\T A\otimes_AT$ is a tilting complex in $D^b(\T A)$ with endomorphism ring $\T B$.
\item Let $X$ be a two-sided tilting complex in $D^b(\T A\otimes_k\T B)$ so that $\T A\otimes_AT\simeq X$ in $D^b(\T A)$. Then $A\otimes_{\T A}X$ and $X\otimes_{\T B}B$ are two-sided tilting complexes in $D^b(A\otimes_kB)$. Moreover, $T\simeq A\otimes_{\T A}X$.
\end{enumerate}
\end{Prop}

Proof. Let $T$ be a tilting complex in $D^b(A)$ with endomorphism ring $B$.
Then by
Rickard's theorem \cite[Corollary 5.4]{Rickard} the complex $\T A\otimes_AT$
is a tilting complex in
$D^b(\T A)$ with endomorphism ring $\T B$. By Keller's theorem
\cite{Kellercoderivations}
there is a two-sided tilting complex $X$ in $D^b(\T A\otimes_k\T B^{op})$
so that $_{\T A}|X\simeq \T A\otimes_AT$.

Now,
\begin{eqnarray*}
Hom_{D^b(A)}(A\otimes_{\T A}X,A\otimes_{\T A}X)&=&
Hom_{D^b(A)}(A\otimes_{\T A}\ _{\T A}|X,A\otimes_{\T A}\ _{\T A}|X)\\
&\simeq &Hom_{D^b(A)}(A\otimes_{\T A}\T A\otimes_AT,A\otimes_{\T A}\T A\otimes_AT)\\
&\simeq &Hom_{D^b(A)}(T,T)\\
&\simeq &B
\end{eqnarray*}
as rings.
But, we know that $\T B$ acts on $A\otimes_{\T A}X$ by multiplication on the right.
Hence, we get a ring homomorphism
$$\T B\lra Hom_{D^b(A)}(A\otimes_{\T A}X,A\otimes_{\T A}X).$$
Since we have seen that the endomorphism ring of $A\otimes_{\T A}X$ is
isomorphic to $B$, the mapping
$$\T B\lra Hom_{D^b(A)}(A\otimes_{\T A}X,A\otimes_{\T A}X)$$
factorises through the canonical projection
$\T B\lra B$. Hence, the action of $\T B$ on $A\otimes_{\T A}X$
has $B^*$ in the kernel.

Now, the $\T B\otimes_k \T B^{op}$-module structure of
$Hom_{D^b(A)}(A\otimes_{\T A}X,A\otimes_{\T A}X)$ is the following.
The action of $\T B$ on the first argument gives the $\T B$-action from the left on
$Hom_{D^b(A)}(A\otimes_{\T A}X,A\otimes_{\T A}X)$ and the action of $\T B$ on the right
comes from the action of $\T B$ on the second argument. Both have the degree $2$
nilpotent ideal $B^*$ in the kernel and so, the natural action of
$\T B\otimes_k\T B^{op}$ is actually an action of $B\otimes_kB^{op}$.

Therefore,
$$Hom_{D^b(A)}(A\otimes_{\T A}X,A\otimes_{\T A}X)\simeq B$$
as $B\otimes_kB^{op}$-modules.

Hence $A\otimes_{\T A}X$ is invertible from the left in the sense that
$$Hom_{D^b(A)}(A\otimes_{\T A}X,A\otimes_{\T A}X)\simeq
Hom_{A}(A\otimes_{\T A}X,A)\otimes_A(A\otimes_{\T A}X)\simeq B$$
as $B\otimes_kB^{op}$-modules.

We still need to show that $A\otimes_{\T A}X$ is invertible from the right as well.
Since $\T A$ and $\T B$ are both symmetric, and since $X$ is a two-sided tilting
complex in $D^b(\T A\otimes_k\T B^{op})$, its inverse complex is given by its
$k$-linear dual $Hom_k(X,k)=:\check X$. We claim that $\check X\otimes_{\T A}A$
is a right inverse of $A\otimes_{\T A}X$. Indeed,
by the previous paragraph, we see that $\T B$ acts on the right of $A\otimes_{\T A}X$
via the projection $\T B\lra B$, i.e. the nilpotent ideal $B^*$ is in the kernel
of this action. Analogously, the same holds for $\check X\otimes_{\T A}A$. Hence,
the tensor product over $\T B$ equals the tensor product over $B$ only.
\begin{eqnarray*}
(A\otimes_{\T A}X)\otimes_{B}(\check X\otimes_{\T A}A)
&\simeq&(A\otimes_{\T A}X)\otimes_{\T B}(\check X\otimes_{\T A}A)\\
&\simeq&A\otimes_{\T A}(X\otimes_{\T B}\check X)\otimes_{\T A}A\\
&\simeq&A\otimes_{\T A}\T A\otimes_{\T A}A\\
&\simeq& A
\end{eqnarray*}
as $A\otimes_kA^{op}$-modules.

As a consequence $A\otimes_{\T A}X$ is a two-sided tilting complex with restriction to
the left being $T$. Moreover, it is clear that $\pi_A$ will map $X$ to
$A\otimes_{\T A}X$.
\dickebox

\begin{Theorem}\label{derivedinvarianceofkappahat}
Let $k$ be a perfect field of characteristic $p>0$, let
$A$ and $B$ be finite dimensional  $k$-algebras and suppose
that $D^b(A)\simeq D^b(B)$ as triangulated categories.
Let $F$ be an explicit standard equivalence between $D^b(A)$ and
$D^b(B)$. Then, $F$ induces a sequence of isomorphisms
$HH_m(F):HH_m(A)\lra HH_m(B)$ so that
$$HH_m(F)\circ \hat \kappa_n^{(m);A}=\hat \kappa_n^{(m);B}\circ HH_{p^nm}(F).$$
\end{Theorem}

Proof. Let $T$ be a tilting complex in $D^b(A)$ with endomorphism ring $B$ and let $X$ be a two-sided tilting complex in $D^b(\T A\otimes_k\T B^{op})$ with
$\T A\otimes_AT\simeq X$ in $D^b(\T A)$.
By Proposition~\ref{twosidedtrivialextension} we get $A\otimes_{\T A}X$ is a two-sided tilting complex in $D^b(A\otimes_kB)$.

We now use Theorem~\ref{derivedkappacup}:
$$HH_{m}(F_X)\circ\kappa_n^{(m);\T A}=\kappa_n^{(m);\T B}\circ HH_{p^nm}(F_X).$$
Multiplying with $HH_{p^nm}(\iota_A)$ from the right and with $HH_{m}(\pi_B)$ from the
left gives then
$$HH_{m}(\pi_B)\circ HH_{m}(F_X)\circ\kappa_n^{(m);\T A}\circ HH_{p^nm}(\iota_A)=HH_{m}(\pi_B)\circ\kappa_n^{(m);\T B}\circ HH_{p^nm}(F_X)\circ HH_{p^nm}(\iota_A).$$

We now claim that
$$HH_{p^nm}(F_X)\circ HH_{p^nm}(\iota_A)=HH_{p^nm}(\iota_B)\circ HH_{p^nm}(F_{A\otimes_{\T A}X})$$
and
$$HH_{m}(\pi_B)\circ HH_{m}(F_X)=HH_{m}(F_{A\otimes_{\T A}X})\circ HH_{m}(\pi_A).$$

Recall from Section~\ref{derivedequivonhomology} how $HH_{m}(F_X)$ is defined.
Then, the commutativity relation will be proven as soon as we have that
the diagrams
$$
\begin{array}{ccc}
\left(\check X\otimes_{\T A}\B(\T A)\otimes_{\T A}X\right)
\otimes_{(\T B)^e}\T B
&\lla&
\B\T A\otimes_{(\T A)^e}(X\otimes_{\T B}\check X)\\
\uar\mbox{\scriptsize $HH(\iota_B)$}&&\mbox{\scriptsize $HH(\iota_A)$}\uar\\
\left((\check X\otimes_{\T A}A)\otimes_A\B A\otimes_A(A\otimes_{\T A}X)\right)
\otimes_{B^e}B
&\lla&
\B A\otimes_{A^e}((A\otimes_{\T A}X)\otimes_{B}(\check X\otimes_{\T A}A))
\end{array}
$$
as well as
$$
\begin{array}{ccc}
\left(\check X\otimes_{\T A}\B(\T A)\otimes_{\T A}X\right)
\otimes_{(\T B)^e}\T B
&\lla&
\B\T A\otimes_{(\T A)^e}(X\otimes_{\T B}\check X)\\
\dar\mbox{\scriptsize $HH(\pi_B)$}&&\mbox{\scriptsize $HH(\pi_A)$}\dar\\
\left((\check X\otimes_{\T A}A)\otimes_A\B A\otimes_A(A\otimes_{\T A}X)\right)
\otimes_{B^e}B
&\lla&
\B A\otimes_{A^e}((A\otimes_{\T A}X)\otimes_{B}(\check X\otimes_{\T A}A))
\end{array}
$$
are commutative.

We use the canonical isomorphisms
$$(\check X\otimes_{\T A}A)\otimes_A\B A\otimes_A(A\otimes_{\T A}X)\simeq
\check X\otimes_{\T A}\B A\otimes_{\T A}X$$
and
$$\B A\otimes_{A^e}(A\otimes_{\T A}X\otimes_B\check X\otimes_{\T A}A)
\simeq \B A\otimes_{(\T A)^e}(X\otimes_{\T B}\check X),$$
which proves that we only need to show that the diagrams
$$
\begin{array}{ccc}
\left(\check X\otimes_{\T A}\B(\T A)\otimes_{\T A}X\right)
\otimes_{(\T B)^e}\T B
&\lla&
\B\T A\otimes_{(\T A)^e}(X\otimes_{\T B}\check X)\\
\uar\mbox{\scriptsize $HH(\iota_B)$}&&\mbox{\scriptsize $HH(\iota_A)$}\uar\\
\left(\check X\otimes_{\T A}\B A\otimes_{\T A}X\right)
\otimes_{B^e}B
&\lla&
\B A\otimes_{(\T A)^e}(X\otimes_{\T B}\check X)
\end{array}
$$
and
$$
\begin{array}{ccc}
\left(\check X\otimes_{\T A}\B(\T A)\otimes_{\T A}X\right)
\otimes_{(\T B)^e}\T B
&\lla&
\B\T A\otimes_{(\T A)^e}(X\otimes_{\T B}\check X)\\
\dar\mbox{\scriptsize $HH(\pi_B)$}&&\mbox{\scriptsize $HH(\pi_A)$}\dar\\
\left(\check X\otimes_{\T A}\B A\otimes_{\T A}X\right)
\otimes_{B^e}B
&\lla&
\B A\otimes_{(\T A)^e}(X\otimes_{\T B}\check X)
\end{array}
$$
are commutative.
But this is clear since the unit element of $\T B$ (and of $\T A$ resp.)
is the image of the unit element of $B$ (and of $A$ resp.)
under $\iota$.

Hence the claim is proven.

This implies now
\begin{eqnarray*}
HH_m(F_{A\otimes_{\T A}X})\circ\hat\kappa_n^{(m);A}
&=&HH_{m}(F_{A\otimes_{\T A}X})\circ HH_{m}(\pi_A) \circ\kappa_n^{(m);\T A}\circ HH_{p^nm}(\iota_A)\\
&=&HH_{m}(\pi_B)\circ HH_{m}(F_{X})\circ\kappa_n^{(m);\T A}\circ HH_{p^nm}(\iota_A)\\
&=&HH_{m}(\pi_B)\circ\kappa_n^{(m);\T B}\circ HH_{p^nm}(F_X)\circ HH_{p^nm}(\iota_A)\\
&=&HH_{m}(\pi_B)\circ\kappa_n^{(m);\T B}\circ HH_{p^nm}(\iota_B)\circ HH_{p^nm}(F_{A\otimes_{\T A}X})\\
&=&\hat\kappa_n^{(m);B}\circ HH_{p^nm}(F_{A\otimes_{\T A}X})
\end{eqnarray*}
which proves Theorem~\ref{derivedinvarianceofkappahat}.
\dickebox

\bigskip

\section{Dual numbers as an example}
\label{dualnumbers}

It should be noted that the mappings $\kappa_n^{(m)}$ are not zero in general.
The dual numbers provide an example. In this section we shall show this fact and prove moreover that $\kappa^{(n)}_m=\hat\kappa^{(n)}_m$ in this case.

\medskip

Holm computed the cohomology ring of rings $k[X]/(f(X))$ for all polynomials $f(X)$ and all fields $k$. In particular for $f(X)=X^2$ one obtains a symmetric algebra, the algebra of dual numbers, whose Hochschild homology and cohomology are isomorphic. For fields of characteristic zero Lindenstrauss \cite{Lindenstrauss} computed in general the Hochschild homology of algebras $k[X_1,\dots,X_n]/{\mathfrak m}^m$ for $\mathfrak m$ being the maximal ideal corresponding to the point $(0,0,\dots,0)$ by exhibiting an explicit projective resolution.

We shall reprove parts of these statements since we will need quite detailed information about the
(co-)cycles that represent each element in the Hochschild structure.

Let $A=k[\epsilon]/(\epsilon^2)$ throughout this section and let
$k$ be a field of characteristic $p>0$.

\begin{Lemma}
$\T A\simeq A\otimes_kA$.
\end{Lemma}

Proof.
Clearly, since $A$ is commutative, also $\T A$ is commutative. Moreover, there is a $k$-basis $A=k\cdot 1\oplus k\cdot \epsilon$ and since $A$ is symmetric, the symmetrising form being $\langle x+y\epsilon,x'+y'\epsilon\rangle:=xy'+x'y$, we get a $K$-basis of $\T A=A^*\sdp A$ by
$$\T A=k\cdot (\langle 1,-\rangle,0)\oplus k\cdot (\langle \epsilon,-\rangle,0)
\oplus k\cdot (\langle 0,-\rangle,1)\oplus k\cdot (\langle 0,-\rangle,\epsilon).$$
Let $\delta:=(\langle 1,-\rangle,0)$, $\sigma:=(\langle \epsilon,-\rangle,0)$,
$\varepsilon:=(\langle 0,-\rangle,\epsilon)$ and $1:=(\langle 0,-\rangle,1)$ (remarking that this is still the unit element of the trivial extension algebra).
We verify immediately that $\varepsilon^2=\delta^2=0$. Furthermore,
$\varepsilon\cdot\delta=\delta\cdot\varepsilon=\sigma$.
Hence, $\T A$ is the quotient of the quiver algebra with one vertex and two loops by the relations saying that the two loops are nilpotent of order $2$ and that they commute. This describes exactly the algebra $A\otimes_kA$. \dickebox

\begin{Rem}
If $p=2$, then $A$ is isomorphic to the group algebra of the cyclic group of order $2$.
Then $A\otimes_kA^{op}$ is isomorphic to the group algebra of the Klein four group and
the Hochschild cohomology ring is known. Holm showed \cite[Theorem 3.2.1]{holmhabil}
that in this case
$$HH^*(A\otimes_kA^{op})\simeq A[X,Y]$$
where $X$ and $Y$ are algebraically independent of degree $1$.
\end{Rem}

For any $p>0$ we know how to determine the Hochschild homology in this more general case
as well by the K\"unneth formula. Indeed,
$$HH_*(A\otimes_kA)\simeq HH_*(A)\otimes_kHH_*(A)$$ by
the K\"unneth formula.

We have to study the injection
$$HH_*(A)\lra HH_*(\T A)\simeq HH_*(A\otimes_KA)\simeq HH_*(A)\otimes_KHH_*(A)$$
given by the isomorphism $\T A\simeq A\otimes_KA$, the injection $A\lra \T A$ and the
K\"unneth formula, as well as the projection
$$HH_*(A)\otimes_kHH_*(A)\simeq HH_*(A\otimes_kA)\simeq HH_*(\T A)\lra HH_*(A).$$

The algebra $A$ is symmetric and hence
$$Hom_k(HH_m(A),k)\simeq HH^m(A)$$
by an isomorphism induced by the symmetrising bilinear form. Therefore there
is a non degenerate pairing $$HH^m(A)\times HH_m(A)\lra K$$
as usual.

\begin{Rem}
The algebra $A\otimes A$ is the quotient of the quiver algebra with one vertex and
two loops by the relations saying that the two loops are nilpotent of order $2$ and
that they commute. Hence, we may replace any loop by a non trivial linear combination
of these two loops, completing by another linear combination of the two loops so that
the determinant of the coefficient matrix is non zero. Hence, we may suppose that the
inclusion $A\lra \T A$ is given by the inclusion $A\lra A\otimes A$ defined by
$a\mapsto 1\otimes a$. Indeed, the Hochschild homology computation by
means of the corresponding double complex does not depend on the choice
of a basis.
\end{Rem}

BACH \cite{Bach} and Holm \cite{Holmpolynomial} computed an explicit resolution of bimodules for a monogenic algebra, such as the dual numbers. Abbreviate for simplicity $A=K[\epsilon]/\epsilon^2$ (as usual) and $A^2:=A\otimes_kA$.
Then a free resolution $C$ of $A$ is periodic of period $2$ and is given by
$$C:\;\;\;(A\lla) A^2\stackrel{d_1}{\lla} A^2\stackrel{d_2}{\lla} A^2 \stackrel{d_1}{\lla} A^2\stackrel{d_2}{\lla} A^2\lla \cdots$$
where $d_1$ is multiplication by $1\otimes\epsilon-\epsilon\otimes 1$ and $d_2$ is
multiplication by $1\otimes\epsilon+\epsilon\otimes 1$.

Applying the functor $-\otimes_{A^2}A$ gives a complex
$$A\stackrel{(d_1)_\otimes}{\lla} A\stackrel{(d_2)_\otimes}{\lla} A \stackrel{(d_1)_\otimes}{\lla} A\stackrel{(d_2)_\otimes}{\lla} A\lla \dots$$
with $(d_1)_\otimes=0$ and $(d_2)_\otimes=2\epsilon$.

Now, applying the functor $Hom_{A^2}(-,A)$ gives a complex $Hom_{A^2}(C,A)$
$$\dots\lra A\stackrel{(d_1)_h}{\lra} A\stackrel{(d_2)_h}{\lra} A \stackrel{(d_1)_h}{\lra} A\stackrel{(d_2)_h}{\lra} A$$
where again $(d_1)_h=0$ and $(d_2)_h=2\epsilon$.

The K\"unneth formula  gives that the tensor product of the cohomology is the cohomology of the tensor product
$$H\left(Hom_{A^2}(C,A)\otimes Hom_{A^2}(C,A)\right)=
H\left(Hom_{A^2}(C,A)\right)\otimes H\left(Hom_{A^2}(C,A)\right).$$
We observe that
$$Hom_{A^2}(C,A)\otimes Hom_{A^2}(C,A)\simeq
Hom_{A^2\otimes A^2}(C\otimes C,A\otimes A)$$
and this is the complex computing the Hochschild cohomology of $A\otimes_kA$, and whence for $\T A$.

As usual the Hochschild (co-)homology depends on weather $p=2$ or $p>2$. Our arguments use the explicit structure of the Hochschild (co-)homology and therefore we shall need to treat these two cases separately.

\subsection{The case $p>2$}

Holm \cite{Holmpolynomial} shows that
$$HH^*(A)\simeq A[U,Z]/(Z\epsilon, U\epsilon, U^2)$$
for an element $U$ in degree $1$ and an element $Z$ in degree $2$. Hence
$HH^m(A)$ is one-dimensional, generated by $Z^n$ if $m=2n$ or $UZ^n$ if $m=2n+1$,
for $m>0$ and isomorphic to $A$ in degree $0$.

We may therefore choose for each $n\in{\mathbb N}\setminus \{0\}$ an element $z_{n}$ in $HH_{2n}(A)$ which corresponds to $Z^n$ in $HH^*(A)$ under the above isomorphism and we get $HH_{2n}(A)=kz_n$. Moreover, choose an element $uz_n\in HH_{2n+1}(A)$ which corresponds to $UZ^n$ under the above isomorphism. Hence, $HH_{2n+1}(A)=kuz_n$.

Observe that
$$\left(HH_*(A)\otimes HH_*(A)\right)^*\simeq HH_*(A)^*\otimes HH_*(A)^*\simeq
HH^*(A)\otimes HH^*(A)$$
where the first isomorphism is canonical and the second one is induced by the
isomorphism $HH_*(A)^*\simeq HH^*(A)$. Hence {\em using this chain of isomorphisms}
we choose the $k$-basis
$$x_iy_j\mbox{ in bidegree }(2i,2j)\mbox{ of }HH_{2i+2j}(A\otimes A),$$
$$x_ivy_j\mbox{ in bidegree }(2i+1,2j)\mbox{ of }HH_{2i+1+2j}(A\otimes A),$$
$$x_iy_jw\mbox{ in bidegree }(2i,2j+1)\mbox{ of }HH_{2i+2j+1}(A\otimes A),$$
$$x_ivy_jw\mbox{ in bidegree }(2i+1,2j+1)\mbox{ of }HH_{2i+1+2j+1}(A\otimes A),$$
which is the dual basis element to the corresponding monomial in $X^iY^j$, the
monomial $X^iVY^j$, the monomial $X^iY^jW$ or the monomial $X^iVY^jW$ of
$HH^{*+*}(A\otimes A)$ under that isomorphism. The $p$-th power of all these basis
elements is zero, except the element $X^iY^j$, whose $p$-th power is $X^{pi}Y^{pj}$.

The (minimal) projective resolution of $A^2$ used above can most easily be expressed as a double complex $C\otimes C$, as was shown above.

$$
\begin{array}{ccccccccccccccc}
&&A^2&\lla&A^2&\lla&A^2&\lla&A^2&\lla&A^2&\lla\\
&&\uar&&\uar&&\uar&&\uar&&\uar\\
A^2&\lla&A^2\otimes A^2&\lla&A^2\otimes A^2&\lla&A^2\otimes A^2&\lla&A^2\otimes A^2&\lla&A^2\otimes A^2&\lla\\
\uar&&\uar&&\uar&&\uar&&\uar&&\uar\\
A^2&\lla&A^2\otimes A^2&\lla&A^2\otimes A^2&\lla&A^2\otimes A^2&\lla&A^2\otimes A^2&\lla&A^2\otimes A^2&\lla\\
\uar&&\uar&&\uar&&\uar&&\uar&&\uar\\
A^2&\lla&A^2\otimes A^2&\lla&A^2\otimes A^2&\lla&A^2\otimes A^2&\lla&A^2\otimes A^2&\lla&A^2\otimes A^2&\lla\\
\uar&&\uar&&\uar&&\uar&&\uar&&\uar\\
A^2&\lla&A^2\otimes A^2&\lla&A^2\otimes A^2&\lla&A^2\otimes A^2&\lla&A^2\otimes A^2&\lla&A^2\otimes A^2&\lla\\
\uar&&\uar&&\uar&&\uar&&\uar&&\uar
\end{array}
$$

In order to view the image of $z_n\in HH_{2n}(A)$ in $HH_{2n}(\T A)$
we need to first establish the multiplicative structure of $HH^*(A\otimes A)$ in terms
of maps in the double complex.

The element $X$ is the degree $(2,0)$ mapping consisting of the identity on all
homogeneous components except on the borders of this bi-complex, and the element $Y$ is the degree $(0,2)$ mapping identity on all homogeneous components except on the
borders of this bi-complex.

The element $V$ is the degree $(1,0)$ mapping consisting of $\epsilon\otimes 1$ in all degrees, except on the borders of the bi-complex, where it is $0$.
Likewise the element $W$ is the degree $(0,1)$ mapping consisting of $\epsilon\otimes 1$ in all degrees, except on the borders of the bi-complex, where it is $0$.

We claim that the element $X^iY^j$ is represented by the degree $(2i,2j)$ mapping
consisting of the identity on all homogeneous components except on the borders of
this bi-complex. The proof of this statement is an easy induction on $i$ and $j$.
Actually, the composition of a morphism of bi-complexes of degree $(2i,2j)$ of the
given shape by a morphism of complexes of degree $(2,0)$, or $(0,2)$ respectively,
corresponds to a morphism of bi-complexes of degree $(2(i+1),2j)$, or $(2i,2(j+1))$
respectively. This corresponds to the cup product of $X^iY^j$ with $X$, or $Y$
respectively.

The element $X^iVY^j$ is represented by the degree $(2i+1,2j)$ mapping
consisting of the mapping $(\epsilon\otimes id)\otimes (id\otimes id)$
on all homogeneous components except on the borders of this bi-complex. Likewise
the element $X^iY^jW$ is represented by the degree $(2i,2j+1)$ mapping
consisting of the mapping $(id\otimes id)\otimes (\epsilon\otimes id)$
and the element $X^iVY^jW$ is represented by the degree $(2i+1,2j+1)$ mapping
consisting of the mapping $(\epsilon\otimes id)\otimes (\epsilon\otimes id)$.

The cup product on $HH^*(A)$ is similar: $Z$ corresponds to the degree $2$ mapping
being the identity on all homogeneous of $C$, except the degree $0$ and degree $1$
component, where the mapping clearly is $0$. Now, again by an analogous induction
as in the bi-complex case, the element $Z^n$ corresponds to the degree $2n$ mapping
being the identity on all homogeneous of $C$, except the degrees up to $2n-1$ component,
where the mapping clearly is $0$. The element $U$ corresponds to the degree $1$ mapping
$\epsilon\otimes 1$ in all degrees except the degree $0$, where it is $0$. The cup
product with $Z^n$ is just an additional shift in degree, the cup product of $U$
with $U$ is $0$, since $(\epsilon\otimes 1)^2=0$.

We now need to compare this mapping to the dual of the Hochschild homology side.
Hence, we dualise the bi-complex and the cocycle representing $X^iY^j$. This just reverses the arrows, identifying again $A^*$ with $A$. The dual of the cocycle $X^iY^j$ can also be obtained by first applying $Hom_{A^2\otimes A^2}(-,A^2)$
to the double complex and then $K$-dualising the result. We obtain the double complex $D_{HH(A^2)}$.
$$
\begin{array}{ccccccccccccccc}
&&.&&.&&.&&.&&.&\\
&&\uar&&\uar&&\uar&&\uar&&\uar\\
.&\lla&A^2&\stackrel{0}{\lla}&A^2&\stackrel{}{\lla}&A^2&\stackrel{0}{\lla}&A^2&
\stackrel{}{\lla}&A^2&\lla\\
&&\uar 0&&\uar 0&&\uar 0&&\uar 0&&\uar 0\\
.&\lla&A^2&\stackrel{0}{\lla}&A^2&\lla&A^2&\stackrel{0}{\lla}&A^2&\lla&A^2&\lla\\
&&\uar &&\uar&&\uar &&\uar&&\uar \\
.&\lla&A^2&\stackrel{0}{\lla}&A^2&\stackrel{}{\lla}&A^2&\stackrel{0}{\lla}&
A^2&\stackrel{}{\lla}&A^2&\lla\\
&&\uar 0&&\uar 0&&\uar 0&&\uar 0&&\uar 0\\
.&\lla&A^2&\stackrel{0}{\lla}&A^2&\lla&A^2&\stackrel{0}{\lla}&A^2&\lla&A^2&\lla\\
&&\uar&&\uar&&\uar&&\uar&&\uar
\end{array}
$$
which now represents exactly the complex computing Hochschild homology of
$A\otimes A$.
Half of the arrows represent the $0$-morphism, the rest, without a side- or superscript in the diagram, represent the mapping $2\epsilon$. The cycle representing $X^n$ is the $2n$-th mapping identity in the upper line.
The complex computing the Hochschild homology of $A$ is simpler:
$$C_{HH(A)}:\;\;A\stackrel{0}{\lla} A\lla A\stackrel{0}{\lla} A\lla A\stackrel{0}{\lla} \cdots$$
with the same convention on non superscribed arrows.

The injection $A\lra A\otimes A\simeq \T A$ maps $a\mapsto 1\otimes a$ and the  projection $A\otimes A\lra A$ maps $a\otimes b\mapsto ab$. Hence, the projective resolution $P_\bullet$ of $A$ as $A^2$-modules maps to the projective resolution $P_\bullet\otimes P_\bullet$ of $A^2$ as $A^2\otimes A^2$-modules as $x\mapsto 1\otimes x$. The second degree component is $0$, except in degree $0$, where it is constantly $1\otimes 1$. Hence, the Hochschild homology complex
of $A$ injects by the identity into the first line of the Hochschild homology (double-) complex of $A^2$.
In other words, the injection produces the mapping of complexes from
$C_{HH(A)}$ to the first line of the bi-complex $D_{HH(A^2)}$ by the mapping $A\ni a\mapsto 1\otimes a\in A\otimes A$ on the level of each homogeneous component.

Likewise, the projection produces a mapping of complexes in the inverse order, which compose to the identity on $C_{HH(A)}$.
Hence $z_n$ is mapped to $x_n$ for all $n>0$ and $z_nu$ is mapped to $x_nv$ for all $n>0$.

Finally, since $p$ is odd, $p^nm$ is odd if and only if $m$ is odd.
We proved the following

\begin{Prop} Let $k$ be a perfect field of characteristic $p>2$. Then
$\hat\kappa_n^{(m),k[\epsilon]/\epsilon^2}=\kappa_n^{(m),k[\epsilon]/\epsilon^2}$
for all $n$ and $m$.
Moreover, $\kappa_n^{(m)}(z_{p^nm})=z_m$ and $\kappa_n^{(m)}(z_{p^nm-1}u)=0$.
\end{Prop}

\subsection{The case $p=2$}

The case of even characteristic is completely analogous, with a few exceptions.
The generator $Z$ of $HH^*(A)$ is in degree $1$, $Z$ is not annihilated by $\epsilon$,
and the differentials for the homology and the cohomology complex are all $0$.
The Hochschild cohomology ring does only contain
nilpotency coming from the centre. The rest is immediate.

\begin{Prop}
Let $k$ be a perfect field of characteristic $2$. Then
$\hat\kappa_n^{(m),k[\epsilon]/\epsilon^2}=\kappa_n^{(m),k[\epsilon]/\epsilon^2}$
for all $n$ and $m$.
Moreover, $\kappa_n^{(m)}(z_{p^nm})=z_m$ for all $m$ and $n$.
\end{Prop}

{\bf Proof.} The proof is a straightforward analogue
of the proof in the $p>2$-case. \dickebox

\end{document}